\documentclass{article}

\usepackage{amsmath}
\usepackage{amssymb}
\usepackage{amsthm}

\usepackage{times}

\newtheorem{lemma}{Lemma}[section]
\newtheorem{theorem}[lemma]{Theorem}
\newtheorem{definition}[lemma]{Definition}

\author{Henry Towsner}
\title{Convergence of Diagonal Ergodic Averages}
\date{\today, DRAFT}

\begin{document}
\bibliographystyle{apalike}
\maketitle

\section{Introduction}
Tao \cite{Tao} has recently proved the following theorem:
\begin{theorem}[Main Theorem]
  Let $l\geq 1$ be an integer.  Assume $T_1,\ldots,T_l:X\rightarrow X$ are commuting, invertible, measure-preserving transformations of a measure space $(X,\mathcal{B},\mu)$.  Then for any $f_1,\ldots,f_l\in L^\infty(X,\mathcal{B},\mu)$, the averages
$$A_N(f_1,\ldots,f_l):=\frac{1}{N}\sum_{n=1}^Nf_1(T_1^nx)\cdots f_l(T^n_lx)$$
converge in $L^2(X,\mathcal{B},\mu)$.
\end{theorem}

The case $l=1$ is the mean ergodic theorem, and the result can be viewed as a generalization of that theorem.  The $l=2$ case was proven by Conze and Lesigne \cite{CL}, and various special cases for higher $l$ have been shown by Zhang \cite{Zhang}, Frantzikinakis and Kra \cite{FK}, Lesigne \cite{Lesigne}, Host and Kra \cite{HK}, and Ziegler \cite{Ziegler}.

Tao's argument is unusual, in that he uses the Furstenberg correspondence principle, which is traditionally used to obtain combinatorial results via ergodic proofs, in reverse: he takes the ergodic system and produces a sequence of finite structures.  He then proves a related result for these finitary systems and shows that a counterexample in the ergodic setting would give rise to a counterexample in the finite setting.

This paper began as an attempt to translate Tao's argument into a purely infinite one.  The primary obstacle to this, as Tao points out (\cite{TaoBlog1}), is that the finitary setting provides a product structure which isn't present in the infinitary setting.  In order to reproduce it, we have to go by an indirect route, passing through the finitary setting to produce a more structured dynamical system.  The structure needed, however, is not the full measure theoretic product.  What is needed in the finitary setting is a certain disentanglement of the transformations, which amounts to requiring that the underlying set of points be a product of $l$ sets, with the $i$-th transformation acting only the $i$-th coordinate, together with a ``nice'' projection under a certain canonical factor.  We obtain this in the infinitary setting using an argument from nonstandard analysis.

A measure space with this property gives rise to measure spaces on each coordinate, but need not be the product of these spaces: it could contain additional measurable sets which cannot be approximated coordinatewise.  These additional sets turn out to be key to the proof, since they are in some sense ``uniform'': they behave, relative to the commuting transformations, as if they were random.  Perhaps unsurprisingly, the behavior of such sets has turned out to be central to a proof of an infinitary analogue of the hypergraph regularity lemma by Elek and Szegedy \cite{ES}.

Using another nonstandard argument to pass from discrete averages to integrals, we show that the non-random functions can be approximated by certain functions of lower complexity in a certain sense.  Proceeding by induction from low complexity to high complexity, we will be able to prove the theorem, using arguments which are essentially those given in \cite{Tao}, translated to an infinitary setting.  This second nonstandard argument has a Furstenberg-style proof as well; a similar argument, along with other variations, is given in \cite{TowsnerGeneralization}.

We thank Jeremy Avigad for providing many helpful suggestions.  We also thank Terrence Tao for answering a number of questions about his proof, and Tim Austin and John Griesmer for finding significant errors in earlier versions of this proof.

\section{Extensions of Product Spaces}
We wish to reduce convergence of the expression in Theorem 1.1 in arbitrary spaces to convergence in spaces where the transformations have been, in some sense, disentangled.  The useful location turns out to be extensions of product spaces---that is, given an ergodic dynamical system $\mathbb{Y}=(Y,\mathcal{C},\nu,U_1,\ldots,U_l)$, we would like to find a system $\mathbb{X}=(\prod_{i\leq l}X_i,\mathcal{B},\mu,T_1,\ldots, T_l)$ where each $T_i$ acts only on the $i$-th coordinate, but which preserves enough properties of the original system that proving convergence for all $L^\infty(\mathbb{X})$ functions is sufficient to give convergence for all $L^\infty(\mathbb{Y})$ functions.

$\mathbb{X}$ naturally gives rise to a product space, taking $\mathcal{B}_i$ to be the restriction of $\mathcal{B}$ to those sets depending only on the $i$-th coordinate, but we do not require that $\mathcal{B}$ be the product of the $\mathcal{B}_i$; in general, $\mathcal{B}$ may properly extend the product.

Given any such system, there is a maximal factor $\mathbb{X}'=(X',\mathcal{B}',\mu\upharpoonright\mathcal{B}')$ in which all sets are $T_iT_j^{-1}$ invariant for each $i,j\leq l$.  We must either accept poor pointwise behavior, since, for example, this factor does not separate $x$ from $T_iT_j^{-1}x$, or, as will do here, take $X'$ to be a different set.  Formally, we will want the property that if $\gamma$ is the projection of $\prod X_i$ onto $X'$, then for every $i\leq l$ and almost every $x_1,\ldots,x_{i-1},x_{i+1},\ldots,x_l$, the function $x_i\mapsto\gamma(x_1,\ldots,x_l)$ is an isomorphism from $(X_i,\mathcal{B}_i,\mu\upharpoonright\mathcal{B}_i)$ to $\mathbb{X}'$.  This obviously requires that all the $\mathbb{X}_i$ be pairwise isomorphic themselves (and further, that $\mathcal{B}$ be symmetric under any change of coordinates).

This requirement is derived from the behavior exhibited by Tao's finitary setting.  Here the product space is the finite measure space on $\mathbb{Z}_N^l$ and $\mathbb{X}'$ is the finite measure space on $\mathbb{Z}_N$.  The map $\gamma:\mathbb{Z}_N^l\rightarrow\mathbb{Z}_N$ is just the map $x_1,\ldots,x_l\mapsto\sum_i x_i$, which has the property that if we fix $x_i$ for $i\neq k$, the map $x_k\mapsto \sum_{i\neq k}x_i+x_k$ is an isomorphism.

Since $(\prod X_i,\mathcal{B},\mu,T_1,\ldots,T_l)$ is not a true product space, we cannot rely on Fubini's Theorem.  Since we nonetheless wish to integrate over coordinates, we have to rely on the use of certain invariant subsets to produce an analogous property.  If $e\subseteq [1,l]$, we will write $x_e$ for an element of $\prod_{i\in e}X_i$; we also write $\overline e$ for the complement of $e$, and, when $e$ is the singleton $\{i\}$, write $\overline i$ for the complement of $\{i\}$.  Given some $x_e$, if $i\in e$ then $x_i$ denotes the corresponding element of the sequence $x_e$.  Given two such variables, say, $x_e,x_{\overline e}$, will write $x$ for the combination of these two vectors, that is
\[(x_e,x_{\overline e})_i:=\left\{\begin{array}{ll}
(x_e)_i&\text{if }i\in e\\
(x_{\overline{e}})_i&\text{otherwise}
\end{array}\right.\]
Note that this is not simply concatenation.  For instance, if $f$ is a function on $\prod X_i$, we will often write $f(x_{\overline k},z)$ as an abbreviation for $f(x_1,\ldots,x_{k-1},z,x_{k+1},\ldots,x_l)$.

\begin{definition}
  Given a measure space $(\prod_{i\leq l} X_i,\mathcal{B},\mu)$, for $k\leq l$, let $\mathcal{B}_{\overline k}$ be the restriction of $\mathcal{B}$ to those sets of the form $\prod_{i\neq k} B_{i}\times X_k$ where $B_i\subseteq X_i$ (or having symmetric difference of measure $0$ with such a set).
\end{definition}
With respect to $\mathcal{B}_{\overline k}$, we may identify elements of $\prod_{i\leq l} X_i$ with elements of $\prod_{i\neq k}X_i$ by discarding the $k$-th coordinate.

\begin{definition}
  Let $\mathbb{Z}$ a dynamical system.  We say a \emph{measure disintegration} exists for some factor $\pi:\mathbb{Z}\rightarrow\mathbb{Z}'$ if there is a map $z'\mapsto\mu_{z'}$ from $\mathbb{Z}'$ to the space of measures on $\mathbb{Z}$ so that $\mu_{z'}$ is supported on $\pi^{-1}(z')$, the map commutes with the group action (so $\mu_{T_gz'}=T_g\mu_{z'}$ for each $g$ and almost every $z'$), and for any $f\in L^2(\mathbb{Z})$,
\[\int fd\mu=\iint fd\mu_{z'}d\mu'\]
where in particular, the right side is defined.
\end{definition}

This disintegration always exists given certain conditions on $\mathbb{Z}$ \cite{FurstenbergBook}, but in our case it is easier to prove that one exists outright than to arrange for those conditions to hold.  We will be dealing with a dynamical system $\mathbb{X}=(\prod_{i\leq l}X_i,\mathcal{B},\mu,T_1,\ldots,T_l)$ where each $T_i$ acts only on the $i$-th coordinate, and where the measure algebra is a (possibly proper) extension of the product measure $\prod_{i\leq l} (X_i,\mathcal{B}_i,\nu_i)$.  Furthermore, taking $\mathcal{B}_{e}$ to consist of those sets depending only on coordinates in $e$, we will assume the measure disintegration onto $(\prod_{i\neq k}X_i,\mathcal{B}_{\overline k},\mu\upharpoonright\mathcal{B}_{\overline k})$ exists (we denote the corresponding measures in the disintegration $\mu_{k,x_{\overline k}}$).  We want to be able to exchange coordinates, and further, to have an additional, $l+1$-st, coordinate which can ``stand in'' for any of the other coordinates.

This extra coordinate will be a factor $\mathbb{X}'=(X',\mathcal{B}',\mu',T_1',\ldots,T_l')$ of $\mathbb{X}$ such that a measure disintegration exists and the projection $\gamma:\prod_{i\leq l}X_i\rightarrow X'$ is $T_iT_j^{-1}$ invariant almost everywhere (that is, for almost every $x$, $\gamma(x)=\gamma(T_iT_j^{-1}x)$ for each $i,j$).  If we fix all but one coordinate of $\mathbb{X}$, we obtain a function $\gamma_{x_{\overline k}}:X_k\rightarrow X'$ by setting $\gamma_{x_{\overline k}}(x_k):=\gamma(x_{\overline k},x_k)$.

\begin{definition}
  Let $\mathbb{X}$ have the form $(\prod_{i\leq l}X_i,\mathcal{B},\mu,T_1,\ldots,T_l)$, and let $\mathbb{X}'=(X',\mathcal{B}',\mu',T_1',\ldots,T_l')$ be a factor of $\mathbb{X}$ with $\gamma$ the corresponding factor map.  We say $\mathbb{X}'$ \emph{cleanly factors} $\mathbb{X}$ if:
  \begin{itemize}
  \item $\gamma$ is $T_iT_j^{-1}$ invariant for each $i,j$
  \item For each $k\leq l$ and almost every $x_{\overline k}\in\prod_{i\neq k}X_i$, the function $\gamma_{x_{\overline k}}$ is an isomorphism from $(X_k,\mathcal{B},\mu_{k,x_{\overline k}})$ onto $\mathbb{X}'$
  \end{itemize}
\end{definition}
Recall that $\mu_{k,x_{\overline k}}$ is the measure given by the measure disintegration of $(\prod_{i\neq k}X_i,\mathcal{B}_{\overline k},\mu\upharpoonright\mathcal{B}_{\overline k})$ evaluated at $x_{\overline k}$.

It is instructive to consider the finite case, where each $X_i$ is just the discrete measure space on $[1,N]$, as is $X'$, and the function $\gamma(x)$ is simply $\sum_{i\leq l} x_i \mod N$.  The clean factoring property here just asserts that if we fix all but one coordinate (and therefore $y:=\sum_{i\neq k}x_i\mod N$), the map $\gamma_{x_{\overline k}}$, which is given by $\gamma_{x_{\overline k}}(x_k)=y+x\mod N$, is a one-to-one mapping.

\begin{theorem}
  If $\mathbb{Y}=(Y,\mathcal{C},\nu,T_1,\ldots,T_l)$ is an ergodic dynamical system with the $T_i$ commuting, invertible, measure-preserving transformations and $f_1,\ldots,f_l\in L^\infty(\mathbb{Y})$ then there is a dynamical system $\mathbb{X}:=(\prod_{i\leq l}X_i,\mathcal{B},\mu,\tilde T_1,\ldots,\tilde T_l)$ and functions $\tilde f_1,\ldots,\tilde f_l\in L^\infty(\mathbb{X})$ such that:
  \begin{itemize}
  \item For each of the factors $(\prod_{i\neq k}X_i,\mathcal{B}_{\overline k},\mu\upharpoonright\mathcal{B}_{\overline k},\tilde T_1,\ldots,\tilde T_l)$ a measure disintegration exists
  \item An $\mathbb{X}'$ exists which cleanly factors $\mathbb{X}$
  \item For each $i$ there is an $S_i$ such that $\tilde T_i$ has the form $\tilde T_i(x_1,\ldots,x_i,\ldots,x_l)=(x_1,\ldots,S_i x_i,\ldots x_l)$
  \item If $A_N(\tilde f_1,\ldots,\tilde f_l)$ converges in the $L^2$ norm then $A_N(f_1,\ldots,f_l)$ does as well
  \end{itemize}
\label{extProd}
\end{theorem}
Note that in the first $A_N$ above, the transformations in question are the $\tilde T_i$, while in the latter, the transformations are the $T_i$.  The proof depends on arguments from nonstandard analysis and the Loeb measure construction; see, for instance, \cite{Goldblatt} for a reference on these topics.
\begin{proof}
  If $\vec v\in[1,P]^l$, write $T^{\vec v}$ for $T_1^{v_1}\cdots T_l^{v_l}$.  By a multidimensional version of the pointwise ergodic theorem (for instance, the general version of the theorem for amenable groups of transformations \cite{Lindenstrauss}), for any function $g$ and almost every $x$,
\[\int gd\nu=\lim_{P\rightarrow\infty}\frac{1}{P^l}\sum_{\vec v\in[1,P]^l}g(T^{\vec v}x)\]
A point with this property is called \emph{generic} for $g$.  Let $\mathcal{G}$ be the set of polynomial combinations of shifts of the functions $f_i$ with rational coefficients.  Since this is a countable set, we may choose a single point $x_0$ which is generic for every element of $\mathcal{G}$.  For each $g\in\mathcal{G}$, define
\[\hat g(\vec n):=g(T^{\vec n}x_0)\]
Since the $f_i$ are $L^\infty$ functions, we may replace them with functions uniformly bounded by some $M_{f_i}$ only changing them on a set of measure $0$, and we may therefore assume that each $\hat g$ is bounded.

Working in an $\aleph_1$-saturated nonstandard extension, choose some nonstandard $c$.  Using the Loeb measure construction, we may extend the internal counting measure on $[1,c]^l$ to a true external measure $\mu$ on the $\sigma$-algebra generated by the internal subsets of $[1,c]^l$.  The functions $\tilde g:=\hat g^*\upharpoonright[1,c]^l$, the restriction of the nonstandard extension of $\hat g$, are internal, and therefore measurable, and bounded since each $\hat g$ is.

For each $g\in\mathcal{G}$, by the definition of $\mu$
\[\int\tilde gd\mu=st\left(\frac{1}{c^l}\sum_{\vec n\in[1,c]^l}\hat g^*(\vec n)\right)\]
where $st$ is the standard part of a bounded nonstandard real.  Furthermore
\[st\left(\frac{1}{c^l}\sum_{\vec n\in[1,c]^l}\hat g^*(\vec n)\right)=\lim_{P\rightarrow\infty}\frac{1}{P^l}\sum_{\vec v\in[1,P]^l}g(T^{\vec v}x_0)\]
follows by transfer: for any rational $\alpha$ greater than $\lim_{P\rightarrow\infty}\frac{1}{P^l}\sum_{\vec v\in[1,P]^l}g(T^{\vec v}x_0)$ and for large enough $P$, $\alpha$ is greater than the average at $P$, so for all nonstandard $c$, $\alpha$ is greater than the average.  Similarly for $\alpha$ less than the limit.  Putting these together, for any $g\in\mathcal{G}$,
\[\int gd\nu=\int\tilde gd\mu\]

Define
\[\tilde T_i(x_1,\ldots,x_i,\ldots,x_l)=(x_1,\ldots,{(x_i+1)\hspace{-6pt} \mod c},\ldots,x_l)\]
It follows that $\tilde T_i \tilde g=\widetilde{T_i g}$,  and by ordinary properties of limits, $\tilde{\cdot}$ commutes with sums and products.  Therefore in particular,
\[\int\left[A_N(f_1,\ldots,f_l)-A_M(f_1,\ldots,f_l)\right]^2d\nu=\int\left[A_N(\tilde f_1,\ldots,\tilde f_l)-A_M(\tilde f_1,\ldots,\tilde f_l)\right]^2d\mu\]

At each point $x_{\overline k}$ in $(\prod_{i\neq k}[1,c],\mathcal{B}_{\overline k},\mu\upharpoonright\mathcal{B}_{\overline k})$, the Loeb measure construction induces a measure $\mu_{k,x_{\overline k}}$ generated by setting
\[\mu_{k,x_{\overline k}}(B):=st\left(\frac{1}{c}\sum_{n\in [1,c]}\chi_B(x_{\overline k},n)\right)\]
for internal $B$.

Finally, let $\mathbb{X}':=([1,c],\mathcal{B}',\mu')$ be given by the Loeb measure construction on $[1,c]$, and let $\gamma:[1,c]^l\rightarrow[1,c]$ be $\gamma(x_1,\ldots,x_l)=\sum_i x_i\mod c$.  The function $\gamma$ is measurable (since it is internal), and measure-preserving (since it maps exactly $c^{l-1}$ points of $[1,c]^l$ to each point of $[1,c]$).  For each $n\in [1,c]$, we may define
\[\mu'_n(B):=st\left(\frac{1}{c^{l-1}}\sum_{\vec v\in [1,c]^l\mid{\sum v_i=n\mod c}}\chi_B(\vec v)\right)\]
for internal $B$ and extend this to a measure on $\mathcal{B}$ by the Loeb measure construction.

$\mathbb{X}$ is isomorphic to a product of $\mathbb{X}'$ with the Loeb measure on $[1,c]^{l-1}$.  A theorem of Keisler \cite{Keisler1,Keisler2} states that when $U$ and $V$ are hyperfinite sets and $f$ is a measurable function on $Loeb(U\times V)$, the functions $v\mapsto f(u,v)$ are measurable for each $u$ and $\int \int f dv du=\int f d(u\times v)$.  In particular, this means that for any measurable $B$, $\mu(B)=\int \mu_n(B)d\mu'(n)$, so a measure disintegration exists.

For any $k\leq l$ and any $x_1,\ldots,x_{k-1},x_{k+1},\dots,x_l\in \prod_{i\neq k}[1,c]$, $\gamma_{\vec x}$ is a measure-preserving bijection from $[1,c]$ to itself mapping measurable sets to measurable sets, and therefore an isomorphism.
\end{proof}

Using the ergodic decomposition, we may reduce the main theorem to the case where $\mathbb{X}$ is ergodic, and then use Theorem \ref{extProd} to reduce to the following case:
\begin{theorem}
  Let $\mathbb{X}=(\prod_{i\leq l}X_i,\mathcal{B},\mu,T_1,\ldots,T_l)$ be a cleanly factored dynamical system such that each $T_i$ has the form
$$T_i(x_1,\ldots,x_i,\ldots,x_l)=(x_1,\ldots,T'_i x_i,\ldots x_l)$$

Then for any $f_1,\ldots,f_l$ in $L^\infty(\mathbb{X})$, $A_N(f_1,\ldots,f_l)$ converges in the $L^2$ norm.
\end{theorem}

In order to prove this theorem, we need a slightly stronger inductive hypothesis, which is what we will actually prove.
\begin{lemma}
  Let $\mathbb{Y}$ be an arbitrary measure space, and let $\mathbb{X}=(\prod_{i\leq l}X_i,\mathcal{B},\mu,T_1,\ldots,T_l)$ be a cleanly factored dynamical system such that each $T_i$ has the form
$$T_i(x_1,\ldots,x_i,\ldots,x_l)=(x_1,\ldots,T'_i x_i,\ldots x_l)$$

Then for any $f_1,\ldots,f_l$ in $L^\infty(\mathbb{X}\times\mathbb{Y})$, $A_N(f_1,\ldots,f_l)$ converges in the $L^2$ norm.
\label{main}
\end{lemma}

For the remainder of the paper, assume $\mathbb{X}$ has this form and that $\mathbb{X}'$ is the factor witnessing that $\mathbb{X}$ is cleanly factored, and let $\gamma$ be the projection onto this factor.  By restricting to the factor generated by the countably many translations of the functions $f_i$, we may assume $\mathbb{X}$ and $\mathbb{X}'$ are separable.

\section{Diagonal Averages}
Note that the projection $\gamma$ we have constructed is consistent with the transformations $T_i$, in the sense that $\gamma(x)=\gamma(y)$ implies $\gamma(T_ix)=\gamma(T_iy)$.  Furthermore, since $\gamma$ is $T_iT_j^{-1}$-invariant, $\gamma(x)=\gamma(y)$ implies that $\gamma(T_ix)=\gamma(T_jy)$, even if $i\neq j$.
\begin{definition}
  Define $T_{l+1}$ on $\mathbb{X}'$ such that for each $x'\in\mathbb{X}'$, if $\gamma(x)=x'$ then $\gamma(T_ix)=T_{l+1}x'$.
\end{definition}
With the particular construction we have given, this definition makes sense pointwise.  For arbitrary $\mathbb{X}'$ cleanly factoring arbitrary $\mathbb{X}$, this is true only almost everywhere.

We wish to reduce Lemma \ref{main} to the case where $\mathbb{X}$ is ergodic.  In order to apply the usual theorem for the existence of an ergodic decomposition (see \cite{Glasner}), the measure space must be a standard Borel space.  It will be easier to take advantage of the fact that we are working with the $L^2$ norm, and get a weaker ergodic decomposition that suffices for our purposes.  Let $\mathcal{C}$ be the factor consisting of sets which are $T_i$-invariant for each $i$ and fix representations of $E(f\mid\mathcal{C})$ for each $f\in L^2(\mathbb{X})$.  Let $\nu$ be the restriction of $\mu$ to $\mathcal{C}$.  For each point $x\in\prod X_i$, we can define a measure $\mu_x$ by $\int fd\mu_x=E(f\mid\mathcal{C})(x)$ with the property that $\iint fd\mu_x d\nu(x)=\int fd\mu$.  Furthermore, the map $x\mapsto\mu_x$ is $T_i$-invariant for each $i$, since $\mathcal{C}$ is, and $\mu_x$ is ergodic for almost every $x$.  (This argument is the first step of the ordinary proof of the ergodic decomposition, as given in \cite{Glasner}, Theorem 3.42.)

We may carry out the same construction on $\mathbb{X}'$ and observe that this preserves the clean factoring property, so it suffices to prove Lemma \ref{main} in the case where $\mu$ is ergodic.

We wish to extend $\mathbb{X}\times\mathbb{X}'$ to ensure that for each $L^2$ function $f$ on $\mathbb{X}$, the functions $x_{\overline k},x'\mapsto f(x_{\overline k},\gamma^{-1}_{x_{\overline k}}(x'))$ are measurable with integral $\int fd\mu$.  If $\mathbb{X}$ were simply a product $\prod_{i\leq l}\mathbb{X}'$, this would occur automatically.  Since this is not necessarily the case, however, we must copy over all the additional sets by swapping coordinates.

Formally\footnote{John Griesmer suggested this simplified definition of $\mathcal{B}^*$}, take the measure algebra $\mathcal{B}^*$ to make measurable all functions $g(x,x')$ such that $x\mapsto \int g(x,x')d\mu'(x')$ is $L^\infty$, and define
\[\int g(x,x')d\nu=\int g(x,x')d\mu'(x')d\mu\]
Define $\mathbb{X}^*:=(\prod_{i\leq l}X_i\times X',\mathcal{B}^*,\nu)$.

%Since $\mathbb{X}'$ cleanly factors $\mathbb{X}$, it can ``substitute'' for any coordinate.  Consider the map $x\mapsto x_{\overline k},\gamma(x)$ from $\prod X_i$ to $\prod_{i\neq k}X_i\times X'$.  The $\sigma$-algebra consisting of sets of the form $B':=\{(x_{\overline k},x')\mid \gamma^{-1}_{x_{\overline k}}(x')\in B\}$ where $B\in\mathcal{B}$ is isomorphic to $\mathcal{B}$, and we may assign $\lambda_k(B')=\mu(B)$ (the push-forward measure).  Furthermore, taking these algebras to have domain $\prod X_i\times X'$ (where the $i$-th coordinate is ignored), any set belonging to more than one of these algebras, or to one of these algebras in addition to $\mathcal{B}$, must not depend on $x'$, and therefore the measures $\lambda_k$ agree on their common domain.

%Take $\mathcal{B}^*$ to be the $\sigma$-algebra generated by the union of these algebras with $\mathcal{B}$, and extend the measure $\mu$ to a measure $\nu$ on $\mathcal{B}^*$ by taking it to be the relatively independent joining of these algebras over their common factors (see, for instance, \cite{Glasner}).  Define $\mathbb{X}^*:=(\prod_{i\leq l}X_i\times X',\mathcal{B}^*,\nu)$.

Importantly, this retains a measure disintegration onto each coordinate:
\[\int f(x,x',y)d\nu=\iint f(x,x',y)d\nu_{k,(x,x')_{\overline k}}d\nu_k\]
where $\nu_{k,(x,x')_{\overline k}}$ is the pushforward of $\mu_{k,x_{\overline k}}$ under $\gamma$ when $k<l+1$ and $\nu_{l+1,x}(\cdot)$ is $\int \cdot\text{ }\delta_x d\mu$.  %In particular, while this is not a literal product of $\mu$ with $\mu'$, the Fubini theorem
%\[\int f(x,x',y)d\nu=\iint f(x,x',y)d\mu'd\mu\]
%still holds.

\begin{definition}
  By abuse of notation, we take $T_i$, $i\leq l+1$, to be transformations on $\mathbb{X}^*\times\mathbb{Y}$ where $T_i(x,x',y)$ is given by $(T_ix,x',y)$ if $i\leq l$ and $T_{l+1}(x,x',y):=(x,T_{l+1}x',y)$.
\end{definition}

Since we will only refer to the product measures with $\mathbb{Y}$, and to limit the proliferation of measures, henceforth we let $\nu$ denote the measure on $\mathbb{X}^*\times\mathbb{Y}$ and $\mu$ denote the measure on $\mathbb{X}\times\mathbb{Y}$.  We write $\mu_k$ and $\nu_k$ for the restriction of $\mu$ and $\nu$ to the $\sigma$-algebra of $T_k$-invariant sets.

To briefly summarize the construction up to this point, given a measure space $\mathbb{Z}$ and $L^\infty$ functions $f_1,\ldots,f_l$, we have constructed a space $\mathbb{X}^*=(\prod_{i\leq l}X_i\times X',\mathcal{B}^*,\nu,T_1,\ldots,T_{l+1})$ with functions $\tilde f_1,\ldots,\tilde f_l$ such that:
\begin{itemize}
\item Convergence of $A_N(\tilde f_1,\ldots,\tilde f_l)$ implies convergence of $A_N(f_1,\ldots,f_n)$
\item The transformations $T_i$ each act only on the $i$-th coordinate
\item The space $\mathbb{X}^*$ has a measure disintegration onto each coordinate, and onto the space of $T_i$-invariant functions for each $i$
\item There is a function $\gamma:\prod_{i\leq l}X_i\rightarrow X'$ cleanly factoring the space $(\prod_{i\leq l}X_i,\mathcal{B},\mu,T_1,\ldots,T_l)$
\end{itemize}

\begin{definition}
    Let $e\subseteq [1,l+1]$.  We say $f\in L^2(\mathbb{X}^*\times\mathbb{Y})$ is $e$-measurable if it is $T_i$-invariant for each $i\not\in e$.  We define $\mathcal{I}_d:=\{e\subseteq [1,l+1]]\mid |e|=d\}$.  We say $f$ has complexity $d$ if it is a finite sum of functions of the form $\prod_{e\in\mathcal{I}_d}g_e$ where each $g_e$ is $e$-measurable.
\end{definition}

\begin{lemma}
  If $f\in L^2(\mathbb{X}^*\times\mathbb{Y})$ is $e$-measurable for some $e$ with $|e|<l+1$ then $f(x,\gamma(x),y)$ is an $L^2$ function and $||f(x,\gamma(x),y)||_{L^2(\mathbb{X}\times\mathbb{Y})}= ||f||_{L^2(\mathbb{X}^*\times\mathbb{Y})}$.
\end{lemma}
\begin{proof}
  For any $i\not\in e$, 
\begin{align*}
\int [f(x,\gamma(x),y)]^2 d\mu
=&\iint [f(x_{\overline i},x_i,\gamma_{x_{\overline i}}(x_i),y)]^2 d\mu_{i,x_{\overline i}}d\mu_i\\
=&\iint [f(x_{\overline i},\gamma^{-1}_{x_{\overline i}}(x'),x',y)]^2 d\mu'd\mu_i\\
\end{align*}
since $\gamma_{x_i}$ is an isomorphism.  Since $f$ is $T_i$-invariant and $x_i$ is ergodic with respect to $T_i$, this is equal to
\[\iint [f(x_{\overline i},x_i,x',y)]^2 d\mu'd\mu\]
Recall that $x_{\overline i},x_i$ is identical to the vector $x$.  But the measure $\nu$ was constructed so this is precisely
\[\int [f(x,x',y)]^2d\nu\]
\end{proof}

In particular, this means that $f(x,\gamma(x),y)$ is an $L^2$ function when $f\in L^2(\mathbb{X}^*\times\mathbb{Y})$ has complexity $d$ for some $d<l+1$.

\begin{definition}
  If $f\in L^\infty(\mathbb{X}^*\times\mathbb{Y})$ has complexity $d$, define $\Delta_N f\in L^\infty(\mathbb{X}\times\mathbb{Y})$ by
\[\Delta_N f:=\frac{1}{N}\sum_{n=1}^N f(x,T_{l+1}^n \gamma(x),y)\]
\end{definition}

We can reduce the question of the convergence of $A_N$ to the convergence of $\Delta_N$:
\begin{definition}
  If $f\in L^2(\mathbb{X}\times\mathbb{Y})$, define $f^i(x,x',y):=f(x_{\overline i},\gamma^{-1}_{x_{\overline i}}(x'),y)$.
\end{definition}
Note that $f^i(x,T_{l+1}^nx',y)=f(x_{\overline i},\gamma^{-1}_{x_{\overline i}}(T_{l+1}^nx'),y)=f(x_{\overline i},T_i^n\gamma^{-1}_{x_{\overline i}}(x'),y)$.

\begin{lemma}
  Let $f_1,\ldots,f_l$ be given.  $A_N(f_1,\ldots,f_l)$ converges in the $L^2$ norm iff $\Delta_N\prod_{i\in\{1,\ldots,l\}} f^i_i$ converges in the $L^2$ norm.
\end{lemma}
\begin{proof}
\begin{align*}
\Delta_N\prod f^i_i(x,y)
=&\frac{1}{N}\sum_{n=1}^N\prod_i f^i_i(x,T^n_{l+1}\gamma(x),y)\\
=&\frac{1}{N}\sum_{n=1}^N\prod_i f_i(x_{\overline i},\gamma^{-1}_{x_{\overline i}}(T^n_{l+1}\gamma(x)),y)\\
=&\frac{1}{N}\sum_{n=1}^N\prod_i f_i(x_{\overline i},T^n_i\gamma^{-1}_{x_{\overline i}}(\gamma(x)),y)\\
=&\frac{1}{N}\sum_{n=1}^N\prod_i f_i(x_{\overline i},T^n_ix_i,y)\\
=&A_N(f_1,\ldots,f_l)(x,y)
\end{align*}
\end{proof}

Each $f^i_i$ is $[1,l+1]\setminus\{i\}$-measurable, so to prove the main theorem, it suffices to prove convergence of $\Delta_N g$ for functions of complexity $d<l+1$.

While $\Delta_N f$ was defined as a function in $L^\infty(\mathbb{X}\times \mathbb{Y})$, we will sometimes view it as the function in $L^\infty(\mathbb{X}^*\times \mathbb{Y})$ where $x'$ is a dummy variable.

\begin{lemma}
  If $g$ and $f$ are $L^\infty(\mathbb{X}^*\times\mathbb{Y})$ functions with complexity $d<l+1$ and $g$ is $T_{l+1}$-invariant then $\Delta_N gf=g\Delta_N f$.
\end{lemma}
\begin{proof}
  Immediate from the definition.
\end{proof}

\begin{lemma}
  Suppose $g$ has complexity $1$.  Then $\Delta_N g$ converges in the $L^2$ norm.
\end{lemma}
\begin{proof}
  If for almost every $y\in Y$, we have convergence for $x\mapsto g(x,y)$ then we may apply the dominated convergence theorem to obtain convergence over $\mathbb{X}^*\times\mathbb{Y}$.  Since $\Delta_N$ distributes over sums, we may further assume that $g$ has the form $\prod_i g_i$ where each $g_i$ is $\{i\}$-measurable.  Then $\Delta_N g=\prod_{i\neq l+1}g_i\Delta_N g_{l+1}$, and by the previous lemma, it suffices to show that $\Delta_N g_{l+1}$ converges.  But this follows immediately from the mean ergodic theorem.
\end{proof}

Because the inductive step generalizes the proof of the ordinary mean ergodic theorem, it is instructive to consider the form of that proof.  The key step is proving that the function $g_{l+1}$ can be partitioned into two components; these components are usually described as an invariant component $g_\bot$ and a component $g_\top$ in the limit of functions of the form $u-T_{l+1}u$.  Unfortunately, this characterization of the second set does not generalize.  There is an alternative characterization, namely that $g_\top$ has the property that $||\Delta_N g_\top||$ converges to $0$.  This turns out to be harder to work with (and, in particular, this characterization does not seem to give a pointwise version of the theorem), but it can be extended to a higher complexity versions.

We will argue as follows: take a function of complexity $d$ in the form $\prod g_e$ with each $g_e$ $e$-measurable, and argue that each $g_e$ can be written in the form $g_{e,\bot}+g_{e,\top}$, where $g_{e,\top}$ is suitably random, so that $||\Delta_N g_{e,\top}\prod h_{e'}||\rightarrow 0$, while $g_{e,\bot}$ is essentially of complexity $d-1$.  If we observe that constant functions have complexity $0$, the usual proof of the mean ergodic theorem has the same form.

\section{From Averages to Integrals}
We need a way to pass from discrete limits to an integral in order to apply the inductive hypothesis.

\begin{lemma}
  Let $\mathbb{X}=(X,\mathcal{B},\mu)$ be a separable measure space and let $b$ be a real number. For $s\leq k$, let $\mathbb{X}_s$ be a factor of $\mathbb{X}$ and $\{b_{m,m',s}\}_{m\leq m'\in\mathbb{N}}$ be a sequence of $L^\infty(\mathbb{X}_s)$ functions bounded (in the $L^\infty$ norm) by $b$.  Let $\{m_t\}_{t\in\mathbb{N}}$ be a sequence such that
$$\frac{1}{m_t}\sum_{i=1}^{m_t}\prod_{s\leq k}b_{i,m_t,s}$$
converges weakly to $f$.  Then there is a space $\mathbb{Y}=(Y,\mathcal{D},\sigma)$ and functions $\tilde b_s\in L^\infty(\mathbb{X}_s\times \mathbb{Y})$ such that $f(x)=\int \prod \tilde b_s(x,y)d\sigma$ for almost every $x$.
\label{floeb}
\end{lemma}
\begin{proof}
  Consider an $\aleph_1$-saturated nonstandard extension of a universe containing $\mathbb{X}$ and the sequences $\{m_t\}$ and $\{b_{m,m_t,s}\}$.  %For convenience, we assume that the extension is obtained by an ultrapower construction.
  Then for each $s\leq k$, there is a nonstandard extension of the sequence $\{b_{m,m_t,s}\}_{m\leq m_t\in\mathbb{N}}$, which we denote $b^*_{m,m_t,s}$.  Let $M:=m_{t'}$ for some nonstandard $t'$, and let $Y:=[1,M]$.
$Y$ is a hyperfinitely additive measure space (taking the counting measure on $Y$), and so, by the Loeb measure construction, there is an external $\sigma$-additive measure, $\sigma$, extending it; we denote the resulting measure space $Loeb(Y)$.  

The elements $b^*_{m,M,s}$ are in $(L^\infty(\mathbb{X}_s))^*$.  Fix some orthonormal basis $\{f_i\}_{i\in\mathbb{N}}$ for $L^2(X)$ consisting of $L^\infty$ functions.  For convenience, we may assume that for each $s$ and each $i$, either $f_i\in L^2(\mathbb{X}_s)$ or $f_i$ is orthogonal to every element of $L^2(\mathbb{X}_s)$.  Then the nonstandard version of this sequence, $\{f_i\}^*$ is a sequence indexed by $i\in\mathbb{N}^*$ which is an orthonormal basis of $L^2(\mathbb{X}^*)$, and further, for $i\in\mathbb{N}$, the $i$-th element is just $f_i^*$.  Therefore the notation $f_i^*$ for $i\in\mathbb{N}^*$ for elements of this sequence is unambiguous.  Define $b^*_s(x,y)=b^*_{y,M,s}(x)$.  Then $b^*_s=\sum_{i\in\mathbb{N}^*}\alpha_{i,s}(y) f^*_i$ for some $\alpha_{i,s}(y)$.  Furthermore, if $\int b_{y,m_t,s} f_id\mu=0$ for each $y,m_t$ then $\alpha_{i,s}=0$, and this is the case whenever $f_i$ is orthogonal to every element of $L^2(\mathbb{X}_s)$.  Since the function $y\mapsto b^*_{y,M,s}$ is internal, so is each $\alpha_{i,s}$, so $st\circ\alpha_{i,s}$ is measurable with respect to $Loeb(Y)$.  Furthermore, since $f^*_i$ and $f^*_j$ are orthogonal when $i\neq j$,
\[\sum_{i\in\mathbb{N}^*}||\alpha_{i,s}||^2_{L^2}||f_i||^2_{L^2}=||b^*_s||^2_{L^2}\leq ||b^*_s||^2_{L^\infty}=b^2\]
Consider the infinite sum $\sum_{i\in\mathbb{N}}(st\circ \alpha_{i,s}(y)) f_i(x)$.  The norm of each finite initial segment of this sum must also be bounded by $b^2$, and so the infinite sum is convergent.  Therefore we may define $\tilde b_s(x):=\sum_{i\in\mathbb{N}}(st\circ \alpha_{i,s}(y)) f_i(x)$.  Since $\tilde b_s$ is a bounded sum of elements of $L^\infty(\mathbb{X}\times Loeb(Y))$, $\tilde b_s\in L^\infty(\mathbb{X}\times Loeb(Y))$.  Furthermore, since $\alpha_{i,s}=0$ whenever $f_i\not\in L^2(\mathbb{X}_s)$, it follows that $\tilde b_s\in L^\infty(\mathbb{X}\times Loeb(Y))$.

In addition, note that the sum over the tails, $\sum_{i>I}\frac{1}{M}\sum_{y=1}^M|\alpha_{i,s}(y)|$ approaches $0$, so in particular, $st(\sum_{i\in\mathbb{N}^*\setminus\mathbb{N}}\frac{1}{M}\sum_{y=1}^M|\alpha_{i,s}(y)|)=0$.

Let $g\in L^2(\mathbb{X})$.  Then for any $\epsilon$ and all but finitely many $t$,
\[\frac{1}{m_t}\sum_{i=1}^{m_t}\int g\prod_{s\leq k}b_{i,m_t,s}d\mu+\epsilon\geq\liminf_{t\rightarrow\infty}\frac{1}{m_t}\sum_{i=1}^{m_t}\int g\prod_{s\leq k}b_{i,m_t,s}d\mu\]
By transfer, it follows that
\[\frac{1}{M}\sum_{y=1}^{M}\int g^*\prod_{s\leq k}b^*_{y,M,s}d\mu+\epsilon\geq\liminf_{t\rightarrow\infty}\frac{1}{m_t}\sum_{i=1}^{m_t}\int g\prod_{s\leq k}b_{i,m_t,s}d\mu\]
for each $\epsilon>0$, and so
\[\frac{1}{M}\sum_{y=1}^{M}\int g^*\prod_{s\leq k}b^*_{y,M,s}d\mu^*\geq\liminf_{t\rightarrow\infty}\frac{1}{m_t}\sum_{i=1}^{m_t}\int g\prod_{s\leq k}b_{i,m_t,s}d\mu\]

We have
\[st(\frac{1}{M}\sum_{y=1}^{M}\int g^*\prod_{s\leq k}b^*_{y,M,s}d\mu^*)=st(\frac{1}{M}\sum_{y=1}^M\int g^*\prod_{s\leq k}\sum_{i\in\mathbb{N}^*}\alpha_{i,s}(y)f_id\mu^*)\]
We may split the sum into the standard and nonstandard components, giving
\[st(\frac{1}{M}\sum_{y=1}^M\int g^*\prod_{s\leq k}\left[\sum_{i\in\mathbb{N}}\alpha_{i,s}(y)f_i+\sum_{i\in\mathbb{N}^*\setminus\mathbb{N}}\alpha_{i,s}(y)f_i\right]d\mu^*)\]
The product expands into sum of $2^k$ components, which we may distribute out of the integral and across the average and standard part operator.  Then the first component is
\[st(\frac{1}{M}\sum_{y=1}^M\int g^*\prod_{s\leq k}\sum_{i\in\mathbb{N}}\alpha_{i,s}(y)f_id\mu^*)\]
while all other components have the form
\[st(\frac{1}{M}\sum_{y=1}^M\int g^* \left(\sum_{i\in\mathbb{N}^*\setminus\mathbb{N}}\alpha_{i,s_0}(y)f_i\right)\prod_{s\neq s_0}\sum_{i\in D_s}\alpha_{i,s}(y)f_id\mu^*)\]
for some $s_0$, and where each $D_s$ is either $\mathbb{N}$ or $\mathbb{N}^*\setminus\mathbb{N}$.  These components are bounded by
\begin{eqnarray*}
st(\frac{1}{M}\sum_{y=1}^M ||\sum_{i\in\mathbb{N}^*\setminus\mathbb{N}}\alpha_{i,s_0}(y)f_i||_{L^2(\mathbb{X})}\cdot ||g^*\prod_{s\neq s_0}\sum_{i\in D_s}\alpha_{i,s}(y)f_i||_{L^2(\mathbb{X})})\\
\leq st(\frac{1}{M}\sum_{y=1}^M \left(\sum_{i\in\mathbb{N}^*\setminus\mathbb{N}}|\alpha_{i,s_0}(y)|\cdot ||f_i||_{L^2(\mathbb{X})}\right)\cdot ||g^*\prod_{s\neq s_0}\sum_{i\in D_s}\alpha_{i,s}(y)f_i||_{L^2(\mathbb{X})})
\end{eqnarray*}
However, as shown above, $st(\frac{1}{M}\sum_{y=1}^M\sum_{i\in\mathbb{N}^*\setminus\mathbb{N}}\alpha_{i,s_0}^2(y))=0$.  Since $||g^*\prod_{s\neq s_0}\sum_{i\in D_s}\alpha_{i,s}(y)f_i||_{L^2(\mathbb{X})}$ is bounded by $||g^*||_{L^2(\mathbb{X})}\prod_{s\neq s_0}\sum_{i\in D_s}\alpha_{i,s}(y)f_i||_{L^2(\mathbb{X})}\leq b^{k-1}$, it follows that all but the first component must be $0$.  Therefore
\[st(\frac{1}{M}\sum_{y=1}^{M}\int g^*\prod_{s\leq k}b^*_{y,M,s}d\mu^*)=st(\frac{1}{M}\sum_{y=1}^M\int g^*\prod_{s\leq k}\sum_{i\in\mathbb{N}}\alpha_{i,s}(y)f_id\mu^*)\]
But by the definition of $\tilde b_s$, this is equal to
\[st(\frac{1}{M}\sum_{y=1}^M\int g\prod_{s\leq k}\tilde b_sd\mu)\]
It follows that
\[\iint g\prod_{s\leq k}\tilde b_s d\mu d\sigma=st(\frac{1}{M}\sum_{y=1}^{M}\int g(x)\prod_{s\leq k}\tilde b_s(x,y)d\mu)\geq\liminf_{t\rightarrow\infty}\frac{1}{m_t}\sum_{i=1}^{m_t}\int g\prod_{s\leq k}b_{i,m_t,s}d\mu\]

A similar argument applies to the $\limsup$, so we conclude that $\iint g(x)\prod_{s\leq k}\tilde b_s(x,y)d\mu d\sigma=\int g(x)f(x)d\mu$ for any $g$.

%The key property of the Loeb measure is that for any internal function $f$, $\int fd\sigma=st(\frac{1}{M}\sum_{y=1}^M f(y))$.  In particular,
%\begin{eqnarray}
%\iint g(x)\prod \tilde b_s(x,y)d\mu d\sigma
%=&st(\frac{1}{M}\sum_{y=1}^M\int g(x)\prod st(b^*_{y,M,s}(x))d\mu)\\
%=&st(\frac{1}{M}\sum_{y=1}^M\int st(g^*(x)\prod b^*_{y,M,s}(x))d\mu)\\
%=&st(\frac{1}{M}\sum_{y=1}^Mst(\int g^*(x)\prod b^*_{y,M,s}(x)d\mu^*))\\
%=&st(\frac{1}{M}\sum_{y=1}^M\int g^*(x)\prod b^*_{y,M,s}(x)d\mu^*)\\
%=&\int gfd\mu
%\end{eqnarray}
%The last step holds because $M$ was chosen precisely to ensure this.  Since this holds for every $g\in L^2(X)$, it follows that $\int \tilde b_s(x,y)d\sigma=f$.
\end{proof}

\section{The Inductive Step}
We now return to the proof of Theorem \ref{main}.  Let $\mathbb{X}=(\prod_{i\leq l}X_i,\mathcal{B},\mu,T_1,\ldots,T_l)$ cleanly factored by $\mathbb{X}'$ be given, and let $\mathbb{Y}$ be an arbitrary measure space.  Recall that $\mathcal{I}_{n}$ is the set of subsets of $[1,l+1]$ with cardinality $n$.  If $e$ is a subset of $[1,l+1]$, we write $\overline{e}$ for the complement of $e$, that is, $[1,l+1]\setminus e$.  We continue to be concerned with functions belonging to $L^\infty(\mathbb{X}^*)$.

\begin{definition}
  Let $e_0\subseteq[1,l+1]$ contain $l+1$.  $Z_{e_0}$ is the subspace of the $e_0$-measurable functions $g$ such that for every sequence $\langle g_e\rangle_{e\in \mathcal{I}_{|e_0|}\setminus\{e_0\}}$ with each $g_e$ $e$-measurable,
\[||\Delta_N g\prod_e g_e||\rightarrow 0\]
as $N$ goes to $\infty$.

$D_{e_0}$ is the set of $e_0$-measurable functions generated by projections onto the $e_0$-measurable sets of weak limit points of sequences of the form 
\[\frac{1}{N}\sum_{n=1}^N \prod_{i\in {e_0}}b_{i,N}(x_{\overline k},T^n_k \gamma^{-1}_{x_{\overline k}}(x'),x',y)\]
as $N$ goes to infinity, for some $k\not\in e_0$, where each $b_i$ is $[1,l+1]\setminus\{i\}$-measurable.
\end{definition}

\begin{lemma}
  If $g$ is $e_0$-measurable where $l+1\in e_0$, $|e_0|<d+1$, and $g\not\in Z_{e_0}$ then there is an $h\in D_{e_0}$ such that $\int ghd\mu>0$.
\end{lemma}
\begin{proof}
  Let an $e_0$-measurable $g\not\in Z_{e_0}$ be given.  Then there is a sequence $\langle g_{e}\rangle_{e\in\mathcal{I}_{|e_0|}\setminus \{e_0\}}$ where each $g_e$ is $e$-measurable and some $\epsilon>0$ such that $||\Delta_N(g\prod_{e\in\mathcal{I}_{|e_0|},e\neq e_0} g_{e})||>\epsilon$ for infinitely many $N$.  Set $f_N:=\Delta_N(g\prod_{e\in\mathcal{I}_{|e_0|},e\neq e_0} g_{e})$.  For each such $N$, we have
\[\int f_N\Delta_N(g\prod_{e\in\mathcal{I}_{|e_0|},e\neq e_0} g_{e}) d\mu>\epsilon^2\]
This means
\[\int \frac{1}{N}\sum_{n=1}^Nf_N(x,y)g(x,T^n_{l+1}\gamma(x),y)\prod_{e\in\mathcal{I}_{|e_0|},e\neq e_0} g_{e}(x,T^n_{l+1}\gamma(x),y) d\mu>\epsilon^2\]
For each $e\neq e_0$, there is some $i\in e_0\setminus e$, so we may assign to each $g_e$ some $i$ such that $g_e$ is independent of $x_i$ and collect the $g_e$ into terms $b_{i,N}$ (independent on $N$, in fact), each a product of some of the $g_e$, such that $b_i$ is independent of $x_i$.  Since $f_N$ is $[1,l]$-measurable, we may also fold $f_N$ into $b_{l+1,N}$, and we have therefore shown that there exist functions $b_{i,N}$ which are $[1,l+1]\setminus\{i\}$-measurable such that
\[\int \frac{1}{N}\sum_{n=1}^Ng(x,T^n_{l+1}\gamma(x),y)\prod_{i\in e_0}b_{i,N}(x,T^n_{l+1}\gamma(x),y)d\mu>\epsilon^2\]
Choosing some $k\not\in e_0$, and letting $g'(x_{\overline k},x',y):=g(x,x',y)$ for almost any $x_k$, this becomes
\[\int g'(x_{\overline k},x',y)\frac{1}{N}\sum_{n=1}^N\prod_{i\in e_0}b_{i,N}(x_{\overline k},T^n_k\gamma^{-1}_{x_{\overline k}}(x'),x',y)d\nu_k>\epsilon^2\]
for infinitely many $N$.  Choosing a subsequence $S$ of these $N$ such that
\[h':=\lim_{N\in S}\frac{1}{N}\sum_{n=1}^N \prod_i b_{i,N}(x_{\overline k},T^n_k\gamma^{-1}_{x_{\overline k}}(x'),x',y)\]
converges, the projection $h$ of $h'$ onto the $e_0$-measurable sets witnesses the lemma.  (In particular, since $g$ is $e_0$-measurable, $\int ghd\mu=\int gh'd\mu>0$.)

\end{proof}

\begin{lemma}
  Every $e_0$-measurable function $g$ may be written in the form $g_\bot+g_\top$ where $g_\bot\in D_{e_0}$ and $g_\top\in Z_{e_0}$.
\end{lemma}
\begin{proof}
  Consider the projection of $g$ onto $D_{e_0}$.  By the previous lemma, if $g-E(g\mid D_{e_0})$ is not in $Z_{e_0}$ then there is an $h\in D_{e_0}$ such that $\int h(g-E(g\mid D_{e_0}))d\mu>0$; this is a contradiction, so $g-E(g\mid D_{e_0})$ belongs to $Z_{e_0}$.
\end{proof}

We could proceed to show that this decomposition is unique, but this is not necessary for the proof.

\begin{lemma}
  If $g=\prod_{e\in\mathcal{I}_{d+1}}g_e$ and each $g_e\in D_e$ then $\Delta_N g$ converges in the $L^2$ norm.
\label{dconv}
\end{lemma}
\begin{proof}
  For convenience, assume $g$ is in the stricter form $\prod_{e\in\mathcal{I}_{d+1},l+1\in e}g_e$.  This is without loss of generality, since if $h=\prod_{e\in\mathcal{I}_{d+1},l+1\not\in e}g_e$ then we have 
\[\Delta_N h\prod_{e\in\mathcal{I}_{d+1},l+1\in e}g_e=h\Delta_N \prod_{e\in\mathcal{I}_{d+1},l+1\in e}g_e\]

  First, assume each $g_e$ is a basic element of $D_{e}$; that is, there is a function $g'_e$ such that $g_e$ is the projection of $g'_e$ onto $\mathcal{B}_{\overline{e_0}}$ and $g'_e$ is a weak limit of an average of the form
\[\frac{1}{N}\sum_{n=1}^N \prod_i b^e_{i,N}(x_{\overline k},T^n_k\gamma^{-1}_{x_{\overline k}}(x'),x',y)\]
Define $b^e_{i,j,n}:=b^e_{i,j}(x_{\overline k},T^n_k\gamma^{-1}_{x_{\overline k}}(x'),x',y)$.  Then Lemma \ref{floeb} applies, so there exist functions $\tilde b^e_i$ such that
\[g'_e(x_{\overline k},x',y)=\int \prod_i \tilde b^e_i(x_{\overline k},z,x',y)d\sigma\]
Since each $g_e$ is the $e$-measurable projection of this function, we may fold $x_{\overline e,0}$ into $z$, integrating over a larger measure space to give
\[g_e(x_e,x',y)=\int \prod_i \tilde b^e_i(x_e,z',x',y)d\sigma'\]
Since each $g_e$ has this form, and these $\tilde b^e_i$ are $e\setminus\{i\}$-measurable, it follows that $g$ has complexity $d-1$, so the result follows by the inductive hypothesis.

If the $g_e$ are sums of basic elements of $D_e$, the result follows immediately.  If $g_e$ is a limit of such elements, each $g_e$ can be written $g_e^0+g_e^1$ where $g^e_0$ is a finite sum of basic elements of $D_e$ and the norm of $g_e^1$ is small.  Then $\prod g_e=\sum_{E\subseteq \mathcal{I}_d}\prod_{e\in E}g^e_0\prod_{e\not\in E}g^e_1$.  When $E=\mathcal{I}_d$, the result follows from the result for finite sums.  When $E\neq\mathcal{I}_d$, the product contains some $g_e^1$, and sine $g_e^1$ is $e$-measurable, it follows that $||\Delta_N g_e||\leq||g_e||$.  Since the $g_{e'}$ are bounded in the $L^\infty$ norm, $||\Delta_N\prod_e g_e||\leq b\prod_e||g_e||$ for some constant $b$, so $\prod_{e\in E}g^e_0\prod_{e\not\in E}g^e_1$ has small norm if $E\neq\mathcal{I}_d$.
\end{proof}

Using this, it is possible to prove Theorem \ref{main}.  If $g=\prod_{e\in\mathcal{I}_{d+1}} g_e(x,x',y)$ where each $g_e$ is $e$-measurable then it suffices to show convergence at each $y$, since then the dominated convergence theorem implies convergence over the whole space.  When $l+1\not\in e$, we have $\Delta_N g_e f=g_e\Delta_N f$, so it suffices to show that $\Delta_N g$ converges where $g$ has the form
$$\prod_{e\in\mathcal{I}_{d+1},l+1\in e} g_e$$
Then write each $g_e$ as $g_{e,\bot}+g_{e,\top}$.  Expanding the product gives
$$\sum_{E\subseteq \{e\in\mathcal{I}_{d+1}\mid l+1\in e\}}\prod_{e\not\in E}g_{e,\bot}\prod_{e\in E}g_{e,\top}$$
where each $g_{e,\top}$ is in $Z_e$ and each $g_{e,\bot}$ is in $D_e$.  Since $\Delta_N$ distributes over sums, it suffices to show that each summand converges.  When $E$ is non-empty, $\Delta_N\prod_{e\not\in E}g_{e,\bot}\prod_{e\in E}g_{e,\top}$ converges to the $0$ function by the definition of $Z_e$.  When $E$ is empty, Lemma \ref{dconv} applies.

\bibliography{Norm}
\end{document}